\documentclass[11pt]{amsart}
\usepackage{mathrsfs}
\usepackage{amsmath}
\usepackage{amscd}
\usepackage{amsmath}

\usepackage{mathrsfs}
\usepackage{amsfonts}
\usepackage[colorlinks,linkcolor=blue,citecolor=blue, pdfstartview=FitH]{hyperref}

\numberwithin{equation}{section}

\theoremstyle{plain}
\newtheorem{thm}{thm}[section]
\newtheorem{theorem}[thm]{Theorem}
\newtheorem{lemma}[thm]{Lemma}
\newtheorem{corollary}[thm]{Corollary}

\theoremstyle{definition}

\newtheorem{remark}[thm]{Remark}

\newtheorem{definition}[thm]{Definition}

\newtheorem{example}[thm]{Example}

\newtheorem{defn-thm}[thm]{Definition-Theorem}

\newcommand{\dbar}{\overline{\partial}}
\newcommand{\pa}{\partial}
\newcommand{\ov}{\overline}
\newcommand{\xu}{\sqrt{-1}}
\newcommand{\Id}{\text{Id}}

\usepackage{fancyhdr}
\begin{document}
\title{Note on the singular Hermitian metrics with analytic singularities}
\makeatletter
\let\uppercasenonmath\@gobble% disables title uppercase
\let\MakeUppercase\relax% disables author uppercase
\let\scshape\relax% disables section smallcaps
\makeatother

\author{Yongpan Zou}

\address{Yongpan Zou, Graduate School of Mathematical Science, The University of Tokyo, 3-8-1 Komaba, Meguro-Ku, Tokyo 153-8914, Japan} \email{zouyongpan@gmail.com}

\maketitle

\begin{abstract}
we study the sheaf of locally square integrable holomorphic section of vector bundle with semi-positive curved singular Hermitian metric.
We confirm the coherence when its induced determinant metric has analytic singularities.
\end{abstract}

\maketitle
 %\tableofcontents

\section{Introduction}

The multiplier ideal sheaves for singular Hermitian metric on holomorphic line bundle play important role in complex geometry and algebraic geometry. A classical result of Nadel says that the multiplier ideal sheaf associated with a plurisubharmonic function is coherent, see \cite[Proposition 5.7]{Dem12a}. There is a natural higher rank analogue of the multiplier ideal sheaf for singular metric of holomorphic vector bundle.

\begin{definition} [cf. \cite{deC98}]
Let $(E,h)$ be a rank $r$ holomorphic vector bundle with singular metric over the complex manifold $X$. Then we define the locally square integrable sheaf of $\mathcal{O}(E)$ as follows, we denote it by $\mathcal{E}(E,h)$,
$$ \mathcal{E}(E,h)_x := \{s \in \mathcal{O}(E)_x : |s|^2_h ~\text{is locally integrable around} ~x\in X \}.
$$
\end{definition}

The fundamental question is when this sheaf is coherent. This problem has been confirmed in some cases, see \cite{deC98, Hos17, Iwa21, HI20, Ina18, Ina20, Ina22} and so on. Since the usual proof make use of H\"ormander's weighted $L^2$-estimate of $\dbar$-equation, some strong positive condition like (singular) Nakano positivity may be necessary (cf. see Example \ref{examp2}). In \cite[Conjecture 1.1]{Ina22}, T. Inayama ask whether $\mathcal{E}(E,h)$ is coherent under the assumption that $h$ is only singular Griffiths semi-positive, i.e., semi-positive curved. He gives a partial answer to this question by assuming the unbounded locus of $\det h$ is discrete.

The determinant line bundle $\det E := \wedge^r E$ of the vector bundle play an essential role in this problem, we denote the induced metric on $\det E$ by $\det h$.
In this note, we give another partial answer to Inayama's conjecture. The following theorem is our main result.

\begin{theorem} [=Corollary \ref{maintheorem1}]
Let $(E,h)$ be a holomorphic vector bundle over an $n$-dimensional complex manifold with a singular Griffiths semi-positive Hermitian metric $h$. If the weight of induced metric $\det h$ on the determinant line bundle $\det E$ have analytic singularities, then $\mathcal{E}(E,h)$ is coherent.
\end{theorem}

\begin{example}
We introduce one important example of G. Hosono. For holomorphic vector bundle $E$, If there exist some global sections $s_1, s_2, \cdots, s_N \in H^0(X, E)$ generically generate $E$, the following morphism of bundle
$$
p : X\times \mathbb{C}^N \rightarrow E
$$
is surjective on a Zariski open subset of X. The induced quotient metric $h$ of $E$ from the standard metric on $\mathbb{C}^N$ is semi-positive curved. As the calculation in \cite[Lemma 4.3]{Hos17}, the determinant metric has analytic singularities. Thus $\mathcal{E}(E,h)$ is coherent by the above Theorem, for details see \cite[Theorem 1.1]{Hos17}.
\end{example}

\begin{example} \label{examp2}
We give two important examples of the sheaf of locally square integrable section. The first come from M. Iwai's paper \cite[Theorem 1.2]{Iwa21}. Let $(E, h)$ be a holomorphic vector bundle on $X$ with a singular Hermitian metric. Iwai proves that $\mathcal{E}(E,h)$ is coherent under the following three conditions:
\begin{enumerate}
\item \label{11} There exists a proper analytic subset $Z$ such that $h$ is smooth on $X\backslash Z$.
\item \label{22} the metric $he^{-\zeta}$ is a positively curved singular Hermitian metric on $E$ for some continuous function $\zeta$ on $X$.
\item \label{33} There exists a real number $C$ such that $\xu\Theta_{E,h}- C\omega\otimes \Id_E\geq0$ on $X\backslash Z$ in the sense of Nakano.
\end{enumerate}
According to Lemma $2.3$ in \cite{Iwa21}, we can replace the second condition with the so-called $L^2$-adapted condition. Recall a local holomorphic frame $v_1,\cdots,v_r$ of vector bundle $(E,h)$ is $L^2$-adapted if for every measurable function $f_1,\cdots,f_r$, the section $\sum f_i v_i$ is locally square integrable if and only if $f_i v_i$ is locally square integrable for every $i=1,\cdots,r$.

In \cite[Corollary B]{SY20}, Schnell and Yang prove the similar result. Let $X$ be a complex manifold and let $D$ be an arbitrary divisor on $X$. Let $\mathcal{V}$ be a polarized variation of rational Hodge structures (VHS) on $X\backslash D$. M. Saito's mixed Hodge module theory shows that $\mathcal{V}$ uniquely corresponds to a polarizable Hodge module $\mathcal{M}$ on $X$ with strict support.  Let $E$ be the lowest nonzero piece in the Hodge filtration of $\mathcal{V}$ and let $(\mathcal{M},F_{\bullet}\mathcal{M})$ be the filtered $D_{X}$-module underlying $\mathcal{M}$. Saito also shows that $E$ extends to the lowest nonzero piece of $F_{\bullet}\mathcal{M}$, which is a torsion-free sheaf on $X$. The Hodge metric $h$ on $E$ extends to a singular Hermitian metric. Let $j:X\setminus D \hookrightarrow X$ be the open embedding. Schnell--Yang proved an interesting fact as followed. Let $\mathcal{F}$ be the subsheaf of $j_{\ast}E$ consisting of sections of $E$
which are locally $L^2$ near $D$ with respect to the Hodge metric on $E$ and the standard Lebesgue measure, then $\mathcal{F}$ is coherent.

It is worth noting that Iwai's result implies Schnell--Yang's result. Indeed, since $E$ is the lowest nonzero piece in the Hodge filtration, by W. Schmid's curvature calculation \cite[Lemma 7.18]{Sch73} it is Nakano semi-positive. On the other hand, according to \cite[Proposition 2.6]{ShZ21}, the vector bundle $j_{\ast}E$ has one $L^2$-adapted holomorphic frame with respect to the extended singular Hermitian metric. Therefore all three conditions above are satisfied.
\end{example}

\noindent\textbf{Acknowledgement}: The author would like to thank his advisor Professor Shigeharu Takayama for guidance, and Professor Junyan Cao for his enlightening question. The author also thanks professor Takahiro Inayama for the helpful discussion.

%\noindent\textbf{Acknowledgement}: The author would like to thank his advisor Professor Shigeharu Takayama for guidance and warm encouragement.

%, and Professor Sheng Rao for his constant support. The author is grateful to the University of Tokyo for Special Scholarship for International Students (Todai Fellowship).

\section{Preliminaries}
In this section, we introduce some basic definitions and results in complex geometry. Unless otherwise mentioned, $X$ denotes a complex manifold of dimension $n$, $\Delta$ denotes a polydisc in $\mathbb{C}^n$. The basic reference is \cite{Dem12b}.
Firstly, let us recall the concept of Chern connection and curvature form of vector bundle.
Let $(E, h)$ be a holomorphic vector bundle on $X$. Corresponding to this metric $h$, there exists the unique Chern connection $D=D_{(E,h)}$, which can be split in a unique way as a sum of a $(1,0)$ and a $(0,1)$ connection, i.e., $D=D'_{(E,h)} + D''_{(E,h)}$.
Furthermore, the $(0,1)$ part of the Chern connection $D''_{(E,h)} =\dbar$. The curvature form is defined to be $\Theta_{E,h} := D^2_{(E,h)}$.
On a coordinate patch $\Omega \subset X$ with complex coordinate $(z_1,\cdots,z_n)$, denote by $(e_1,\cdots,e_r)$ an orthonormal frame of vector bundle $E$ with rank $r$. Set
$$
\xu \Theta_{E,h} = \xu \sum_{1\leq j,k\leq n, 1\leq \lambda,\mu \leq r} c_{jk\lambda\mu} dz_j\wedge d\ov{z}_k \otimes e^{\ast}_{\lambda} \otimes e_{\mu} , \quad c_{jk\mu\lambda}=\ov{c}_{jk\lambda\mu}.
$$
Corresponding to $\xu\Theta_{E,h}$, there is a Hermitian form $\theta_{E,h}$ on $TX\otimes E$ defined by
$$
\theta_{E,h}(\phi,\phi) = \sum_{jk\lambda\mu} c_{jk\lambda\mu}(x) \phi_{j\lambda}\ov{\phi}_{k\mu}, \quad \phi \in T_{x}X\otimes E_x.
$$

\begin{definition}
A holomorphic vector bundle $(E,h)$ is said to be
\begin{enumerate}
    \item \emph{Nakano positive} (resp. \emph{Nakano semi-positive}) if for every nonzero tensor $\phi \in TX\otimes E$, we have
$$
 \theta_{E,h}(\phi,\phi) >0 \quad(\text{resp.} \geq 0).
$$
  \item \emph{Griffiths positive} (resp. \emph{Griffiths semi-positive}) if for every nonzero decomposable tensor $\xi\otimes e \in TX\otimes E$, we have
$$
 \theta_{E,h}(\xi\otimes e, \xi\otimes e) >0 \quad(\text{resp.} \geq 0).
$$
\end{enumerate}
It is clear that Nakano positivity implies Griffiths positivity and that both concepts coincide if $r=1$. In the case of line bundle, $E$ is merely said to be positive (resp. semi-positive).
\end{definition}

\noindent The Nakano positivity has a close relationship with the solvability of bundle valued $\dbar$-equation. The next theorem is fundamental.
\begin{theorem} \cite[Theorem 5.1]{Dem12a} \label{dbareq}
Let $X$ be a complete K\" ahler manifold with a K\"ahler metric $\omega$ which is not necessarily complete. Let $(E,h)$ be a Hermitian vector bundle of rank $r$ over $X$, and assume that the curvature operator $B:=[i\Theta_{E,h},\Lambda_\omega]$ is semi-positive definite everywhere on $\wedge^{p,q}T_X^*\otimes E$, for some $q\geq 1$. Then for any form $g\in L^2(X,\wedge^{p,q}T^*_{X}\otimes E)$ satisfying $\bar{\partial}g=0$ and $\int_X\langle B^{-1}g,g\rangle dV_\omega<+\infty$, there exists $f\in L^2(X,\wedge^{p,q-1}T^*_X\otimes E)$ such that $\bar{\partial}f=g$ and
$$\int_X|f|^2dV_\omega\leq \int_X\langle B^{-1}g,g\rangle dV_\omega.$$
\end{theorem}
In general, the Griffiths positive condition is not enough to solve the $\dbar$-equation of vector bundle, but Demailly--Skoda have the following interesting result.

\begin{theorem} [\cite{DeS79}] \label{De-Sk}
Let $(E,h)$ be a holomorphic vector bundle with smooth Hermitian metric $h$, if $(E,h)$ is Griffiths semi-positive, then $(E\otimes \det E, h\otimes \det h)$ is Nakano semi-positive.
\end{theorem}

%For the definition of singular Hermitian metrics, see \cite[Section 3]{BP08} or \cite[Definition 1.1]{Rau15}.

Next, we come to the singular category, we first introduce positivity notions for singular Hermitian metrics. Let $H_r$ be the space of semi-positive, possibly unbounded Hermitian forms on $\mathbb{C}^r$. A singular Hermitian metric $h$ on vector bundle $E$ is a measurable map from $X$ to $H_r$ such that $h(x)$ is finite and positive definite almost everywhere. In particular we have $0< \det h < +\infty$ almost everywhere.

\begin{definition}
Let $(E,h)$ be a singular Hermitian metric on $X$, then $(E,h)$ is said to be:
\begin{enumerate}
		\item Griffiths semi-negative (or semi-negative curved)if $\log |s|^2_h$ (or $|s|^2_h$) is plurisubharmonic (psh, for short) for any local holomorphic section $s$ of $E$.
		\item Griffiths semi-positive (or semi-positive curved)if the dual metric $h^\star$ on $E^\star$ is Griffiths semi-negative.
	\end{enumerate}
\end{definition}

The plurisubharmonic function is one of the essential concepts in complex geometry. A quasi-plurisubharmonic (quasi-psh, for short) function is a function $v$ which is locally equal to the sum of a psh function and of a smooth function. The following regularization lemma is very useful.

%\Begin{Definition} [Psh Function And Quasi-Psh]
%A Function $U: \Omega \Rightarrow [-\Infty, \Infty)$ Defined On A Open Subset $\Omega \In \Mathbb{C}^N$ Is Called Plurisubharmonic (Psh, For Short) If
%\Begin{Enumerate}
%\Item $U$ Is Upper Semi-Continuous;
%\Item For Every Complex Line $Q \Subset \Mathbb{C}^N$, $U|_{\Omega \Cap Q}$ Is Subharmonic On $\Omega \Cap Q$.
%\End{Enumerate}
%
%\End{Definition}

\begin{lemma} [\cite{BP08, Rau15}] \label{regular}
\begin{enumerate}
\item Suppose $X$ is a polydisc in $\mathbb{C}^n$, and suppose $h$
is a singular Hermitian metric on $E$ which is semi-negative (resp. semi-positive)
curved. Then, on any smaller polydisc, there exists a sequence of smooth
Hermitian metrics $\{h_\nu\}$ decreasing (resp. increasing) pointwise to $h$ whose
corresponding curvature tensor is Griffiths negative (resp. positive).
\item Suppose a singular Hermitian metric $h$ is semi-positive curved. Then $(-\log \det h) \in L^1_{loc}(X, \mathbb{R})$ and is a psh function.
\end{enumerate}
\end{lemma}

Sometimes it is more natural to consider the sheaf $K_X \otimes \mathcal{E}(E,h)$ instead of square integrable sheaf $\mathcal{E}(E,h)$ is the former has nice functorial property.

\begin{lemma} \cite[Proposition 4.3]{ShZ21} \label{modifi}
Let $\pi: X' \rightarrow X$ be a proper modification of complex manifolds, and $(E,h)$ be the vector bundle on $X$ with possible singular metric, then
$$
\pi_{\ast}(K_{X'}\otimes \mathcal{E}(\pi^{\ast}E, \pi^{\ast}h)) = K_X \otimes \mathcal{E}(E,h).
$$
\end{lemma}

\begin{definition}
A quasi-psh function $u$ will be said to has \emph{analytic singularities} if $u$ can be written locally as
$$
u = \alpha \log (|f_1|^2 + \cdot\cdot\cdot + |f_N|^2) + v
$$
where $\alpha$ is a positive real number, $v$ is a locally bounded function and all $f_i$ are holomorphic function. Moreover, if the coefficient $\alpha$ is a positive rational number, then we say $u$ has \emph{algebraic singularities}.
\end{definition}

%\begin{definition} [Singular metric and curvature current on line bundle]
%Let $(F,h)$ be a holomorphic line bundle on complex manifold $X$ endowed with possible singular Hermitian metric $h$. For any given trivilization $\theta : F|_{\Omega} \simeq \Omega \times \mathbb{C}$ by
%$$ \| \xi \|_h = |\theta(\xi)| e^{-\phi(x)},  \quad x \in \Omega, \xi \in F_x,
%$$
%where $\phi \in L^1_{loc}(\Omega)$ is an arbitrary function, called the weight of the metric. The curvature $\sqrt{-1} \Theta_h(F)$ of $h$ is defined by
%$$  \sqrt{-1} \Theta_h(F) = \sqrt{-1} 2\partial \overline{\partial} \phi.
%$$
%The Levi form $\sqrt{-1}\pa\dbar \phi$ is taken in the sense of distributions and thus the curvature is a $(1,1)$-current but not always a smooth $(1,1)$-form. It is globally defined on $X$ and independent of the choice of trivializations. The curvature $\sqrt{-1} \Theta_h(F)$ of $h$ is said to be positive (resp. semi-positive) if $\sqrt{-1} \Theta_h(F) \textgreater 0$ (resp. $\geq 0$) in the sense of current.
%\end{definition}

%\begin{definition} [Multiplier ideal sheaves]
%Let $\phi$ be a quasi-psh function on a complex manifold $X$, the multiplier ideal sheaf $\mathscr{I}(\phi) \subset \mathcal{O}_X$ is defined by
%$$ \Gamma(U, \mathscr{I}(\phi)) = \{f \in \mathcal{O}_X(U) :\quad |f|^2e^{-2\phi} \in L^1_{loc}(U) \}
%$$
%for every open set $U \subset X$. For a line bundle $(F, h)$, if the local weight of metric $h$ is $\phi$, then we denote the multiplier ideal sheaf interchangeably by $\mathscr{I}(\phi)$ or $\mathscr{I}(h)$.
%\end{definition}

\section{Singular Hermitian metrics with algebraic singularities}

In order to investigate the coherence, we need to solve one special $\dbar$-equation first. The following theorem and its proof are the slight modification of Inayama's theorem $3.1$ in \cite{Ina22}. We present it here for reader's convenience. Some simplified notations: $\det^m E:= (\det E)^{\otimes m}, \det^m h:= (\det h)^m$.

\begin{theorem} \label{debarrq1}
Let $(X, \omega)$ be a Stein manifold with a standard K\"ahler metric $\omega$, for any positive constant $k$, we can find a smooth strictly psh function $\phi$ with $\xu \pa\dbar \phi \geq k \omega$. Let $(E,h)$ be the trivial holomorphic vector bundle with a Griffiths semi-positive singular Hermitian metric $h$. There exists a natural induced metric $h\otimes \det^m h$ on the vector bundle $E\otimes \det^m E$ for any $m \in \mathbb{N}$. For any $\dbar$-closed $E\otimes \det^m E$ valued $(n,q)$-form $u$ with finite $L^2$ norm with respect to $h\otimes \det^m h$, there exists a $E\otimes \det^m E$ valued $(n,q-1)$-form $\alpha$ such that $\dbar \alpha = u$ and
$$
\int_X |\alpha|^2_{\omega, h\otimes \det^m h} e^{-\phi} dV_{\omega} \leq \frac{1}{qC} \int_X |u|^2_{\omega, h\otimes \det^m h} e^{-\phi} dV_{\omega}
$$
for some constant $C$.
\end{theorem}

\begin{proof}
By the assumption of Stein, we may regard $X$ as a submanifold of $\mathbb{C}^N$ for some positive integer $N$. We denote by $i: X \hookrightarrow \mathbb{C}^N$ the inclusion of $X$ and according to Siu's result in \cite{Siu76}, there exist an open neighborhood $M$ of $X$ in $\mathbb{C}^N$ and a holomorphic retraction $p: M \rightarrow X$ such that $p \circ i = \Id_X$. Now $(p^{\ast}E, p^{\ast}h)$ is the trivial vector bundle with Griffiths semi-positive metric $p^{\ast}h$ on $M$. From the Lemma \ref{regular}, one can obtain a sequence of smooth Hermitian metrics $\{g_\nu\}$ with Griffiths semi-positive curvature increasing to $p^{\ast}h$ on any relative compact subset in $M$. Taking an exhaustion $\{X_j\}_{j=1}^{\infty}$ of $X$, where each $X_j$ is a relative compact Stein sub-domain in $X$, satisfying each $X_j$ is the relatively compact subset of $X_{j+1}$ and  $\bigcup X_j = X$. Set $\{h_\nu := i^{\ast}g_\nu\}$ and so $\{h_\nu\}$ is an approximation sequence with Griffiths semi-positive curvature increasing to $h$ on any relatively compact subset of $X$. Since each $h_\nu$ is Griffiths semi-positive, due to Demailly-Skoda's theorem \ref{De-Sk}, $h_\nu \otimes \det h_\nu$ is Nakano semi-positive. The curvature of $(E\otimes \det^m E, h_\nu\otimes \det^m h_\nu \cdot e^{-\phi})$ can be calculated as
\begin{align*}
&\xu \Theta_{h_\nu\otimes \det^m h_\nu\cdot e^{-\phi}} \\
=& \xu \Theta_{h_\nu\otimes \det^m h_\nu} + \xu \pa\dbar \phi \otimes \Id_{E\otimes \det^m E} \\
=& \xu \Theta_{h_\nu\otimes \det h_\nu} \otimes \Id_{\det^{m-1} E} + (m-1)\xu \pa\dbar (-\log \det h_\nu) \otimes \Id_{E\otimes \det E} \\
&+ \xu \pa\dbar \phi \otimes \Id_{E\otimes \det^m E} \\
\geq & ~ C \omega\otimes \Id_{E\otimes \det^m E}.
\end{align*}
The last inequality is in the sense of Nakano and $C$ is independent to $\nu$. Thus for any $E\otimes \det^m E$ valued $(n,q)$-form $u$ with finite norm, we have
$$
\langle  [\xu \Theta_{h_\nu\otimes \det^m h_\nu \cdot e^{-\phi}}, \Lambda_{\omega}]u, u \rangle \geq qC |u|^2_{h_\nu\otimes \det^m h_\nu \cdot e^{-\phi}}
$$

Now we fix one sub-domain $X_j$. By using Theorem \ref{dbareq}
%for any $\dbar$-closed $E\otimes \det E$-valued $(n,q)$-form with finite norm,
we get a solution $\alpha_\nu$ of the $\dbar$-equation satisfying
\begin{align*}
\int_{X_j}|\alpha_\nu|^2_{\omega, h_\nu\otimes \det^m h_\nu}e^{-\phi} dV_\omega &\leq \frac{1}{qC}\int_{X_j} |u|^2_{\omega, h_\nu\otimes \det^m h_\nu}e^{-\phi} dV_{\omega}\\
&\leq \frac{1}{qC}\int_{X_j} |u|^2_{\omega, h\otimes \det^m h}e^{-\phi} dV_{\omega}\\
&\leq \frac{1}{qC}\int_X |u|^2_{\omega, h\otimes \det^m h}e^{-\phi} dV_{\omega} <+\infty
\end{align*}
for sufficiently large $\nu$. Fix sufficiently large $\nu_0$.
We have that for $\nu \geq \nu_0$
\begin{align*}
\int_{X_j}|\alpha_\nu|^2_{\omega, h_{\nu_0}\otimes \det^m h_{\nu_0}}e^{-\phi} dV_{\omega} &\leq \int_{X_j}|\alpha_\nu|^2_{\omega, h_\nu\otimes \det^m h_\nu}e^{-\phi} dV_\omega\\
&\leq \frac{1}{qC}\int_X |u|^2_{\omega, h\otimes \det^m h}e^{-\phi} dV_{\omega} <+\infty.
\end{align*}
Then $\{ \alpha_\nu\}_{\nu\geq \nu_0}$ forms a bounded sequence on $X_j$ with respect to the norm
$$\int_{X_j}|~\cdot~|^2_{\omega, h_{\nu_0}\otimes \det^m h_{\nu_0}}e^{-\phi} dV_{\omega}.$$
We can get a weakly convergent subsequence $\{\alpha_{\nu_0, k}\}_k$.
Thus, the weak limit $\alpha_j$ satisfies
$$
\int_{X_j}|\alpha_j|^2_{\omega, h_{\nu_0}\otimes \det^m h_{\nu_0}}e^{-\phi} dV_{\omega}\leq \frac{1}{qC}\int_X |u|^2_{\omega, h\otimes \det^m h}e^{-\phi} dV_{\omega} <+\infty.
$$
Next, we fix $\nu_1 > \nu_0$.
Repeating the above argument, we can choose a weakly convergent subsequence $\{\alpha_{\nu_1, k}\}_k \subset \{ \alpha_{\nu_0, k}\}_k$ with respect to $\int_{X_j}|~\cdot~|^2_{\omega, h_{\nu_1}\otimes \det^m h_{\nu_1}}e^{-\phi} dV_{\omega}$.
Then by taking a sequence $\{ \nu_n\}_n$ increasing to $+\infty$ and a diagonal sequence, we obtain a weakly convergent sequence $\{ \alpha_{\nu_k, k}\}_k$ with respect to $\int_{X_j}|~\cdot~|^2_{\omega, h_{\nu_\ell}\otimes \det^m h_{\nu_\ell}}e^{-\phi} dV_{\omega}$ for all $\ell$.
Hence, $\alpha_j$ satisfies
$$
\int_{X_j}|\alpha_j|^2_{\omega, h\otimes \det^m h}e^{-\phi} dV_{\omega}\leq \frac{1}{qC}\int_X |u|^2_{\omega, h\otimes \det^m h}e^{-\phi} dV_{\omega}
$$
thanks to the monotone convergence theorem.
Since the right-hand side of the above inequality is independent of $j$, by using the exactly same argument, we can get an $E\otimes \det^m E$-valued $(n,q-1)$-form $\alpha$ satisfying $\dbar \alpha =u$ and
$$
\int_X |\alpha|^2_{\omega, h\otimes \det^m h}e^{-\phi} dV_{\omega}\leq \frac{1}{qC}\int_X |u|^2_{\omega, h\otimes \det^m h}e^{-\phi} dV_{\omega},
$$
which completes the proof.
%See the argument in \cite[Theorem 1.1 or 1.3]{Ina18} or \cite[Theorem 3.4]{Ina20} for more details.
\end{proof}

\begin{remark} \label{cohe}
From now on, since the coherence is the local property, we assume $X := \Delta := \Delta^n$ be the polydisc in $\mathbb{C}^n$. If possible, we can shrink $\Delta$ to meet our requirements. Let $\omega$ be the standard K\"ahler metric on $\Delta$.
With Theorem \ref{debarrq1} at hand, one can proof the coherence of $\mathcal{E}(E\otimes \det^m E, h\otimes \det^m h)$. To be specific, one can prove that it is generated by the square integrable sections of $E\otimes \det^m E$.
The proof is the same as the procedure in \cite[Proposition 5.7]{Dem12a}.
\end{remark}

\begin{theorem} \label{maintheorem}
Given a holomorphic vector bundle $E$ with Griffiths semi-positive metric $h$. If the weight of Hermitian metric $\det h$ on the determinant line bundle $\det E$ has algebraic singularities, then the sheaf $\mathcal{E}(E, h)$ is coherent. Note that $(E, h)$ always be Griffiths semi-positive.
\end{theorem}

\begin{proof}
We divide the proof into several steps.

\emph{step $1$}.~ The sheaf $\mathcal{E}(E,h)$ is coherent if and only if $K_X \otimes \mathcal{E}(E,h)$ is coherent. We thus consider the coherence of $K_X \otimes \mathcal{E}(E,h)$ since it has the nice functorial property under the proper modification. Indeed, if $\pi: X' \rightarrow X$ be a proper modification, then due to previous Lemma \ref{modifi}, we have
$$
\pi_{\ast} (K_X \otimes \mathcal{E}(\pi^{\ast}E, \pi^{\ast}h)) =K_X \otimes \mathcal{E}(E,h).
$$
According to Grauert's coherence theorem, if $K_X \otimes \mathcal{E}(\pi^{\ast}E, \pi^{\ast}h)$ is coherent, then its direct image $K_X \otimes \mathcal{E}(E, h)$ is coherent too. Note that $(\pi^{\ast}E, \pi^{\ast}h)$ is also singular Griffiths semi-positive, see for example \cite[Lemma 2.3.2]{PT18}.
By the assumption, the weight of induced metric $\det h$ on the determinant line bundle $\det E$ has algebraic singularities, i.e., locally can be written as
$$
-\log \det h = \gamma \log (\sum_j |h_j|^2) + \xi.
$$
Here $\gamma$ is a positive rational number and $\xi$ is bounded function.
Therefore after blow-up or modification, we can assume that
$$
-\log \det h = \sum_j \gamma_j \log |h_j|^2 + \xi.
$$
All $\gamma_j$ are positive rational number, we can write $\gamma_j = \frac{n_j}{m}$, where $n_j, m \in \mathbb{Z}$ and $m\neq 0$. So we have
\begin{align}\label{algebraic}
-m \log \det h = \sum_j n_j \log |h_j|^2 + m \xi
\end{align}
be the local weight of $\det^m h$, the metric of line bundle $\det^m E$. Now locally the metric
\begin{align}
det^m h = \frac{1}{\prod_j |h_j|^{2n_j}} e^{-m\xi}.
\end{align}
Let $H= \prod h_j^{n_j}$ be the corresponded holomorphic function. Locally $|H|^2 \det^m h = e^{-m\psi}$ is positive and bounded function.

\emph{step $2$}.~ As we have said, since the coherence is a local property, we may assume that $X :=\Delta := \Delta^n$ is a small polydisc in $\mathbb{C}^n$. Let $H^0_{2, h\otimes \det^m h}(\Delta, E\otimes \det^m E)$ be the space of holomorphic section $s$ of $E\otimes \det^m E$ on $\Delta$ such that
$$
\int_{\Delta} |s|^2_{h\otimes\det^m h} dV_\omega < +\infty.
$$
We consider the evaluation map
$$
ev: H^0_{2, h\otimes \det^m h}(\Delta, E\otimes det^m E) \times \mathcal{O}_{\Delta} \rightarrow \mathcal{O}(E\otimes det^m E).
$$
We denote by $e \otimes\det^m e$ the image of $ev$. Every coherent $\mathcal{O}_{\Delta}$-sheaf enjoys the Noether property, it is obvious that $e \otimes\det^m e$ is coherent, and by the Remark \ref{cohe} we have the equality
\begin{align} \label{equa1}
( e \otimes det^m e)_x = \mathcal{E}(E\otimes det^m E, h\otimes det^m h)_x.
\end{align}
Similarly, for the line bundle $(\det^m E, \det^m h)$, we denote by $\det^m e$ the image of square integrable sections under the evaluation map, and according to the coherence of $\mathcal{E}(\det^m E, \det^m h)$, we have
\begin{align}
(det^m e)_x = \mathcal{E}(det^m E, det^m h)_x.
\end{align}
For the vector bundle $(E, h)$, we denote by $e$ the image of square integrable sections under the evaluation map as above. We want to show that $e_x = \mathcal{E}(E,h)_x$ for any point $x\in \Delta$. Since $e$ is coherent and $e_x \subseteq \mathcal{E}(E,h)_x$, one just need to check
\begin{align} \label{equal2}
e_x + \mathcal{E}(E,h)_x \cap m_x^k\cdot E_{(x)}  = \mathcal{E}(E,h)_x
\end{align}
for any positive integer $k$, here $E_{(x)}:= \underrightarrow{\lim}_{x\in U}H^0(U, E)$. Indeed, if this is the case, by the Artin--Rees lemma, there exists a positive integer $l$ such that
$$
\mathcal{E}(E,h)_x \cap m_x^k\cdot E_{(x)}  = m_x^{k-l}\cdot(\mathcal{E}(E,h)_x \cap m_x^l\cdot E_{(x)})
$$
holds for any $k> l$. Therefore according to the above equality \eqref{equal2}, one have
$$
\mathcal{E}(E,h)_x = e_x + \mathcal{E}(E,h)_x \cap m_x^k\cdot E_{(x)} \subset e_x + m_x\cdot \mathcal{E}(E,h)_x \subset \mathcal{E}(E,h)_x.
$$
By Nakayama's lemma, one can obtain $e_x = \mathcal{E}(E,h)_x$, which is the desired result.
Now we begin to prove the equality \eqref{equal2}, one inclusion is trivial, we just need to check that $e_x + \mathcal{E}(E,h)_x \cap m_x^k\cdot E_{(x)}  \supseteq \mathcal{E}(E,h)_x$

\emph{step $3$}.  Let $f_x \in \mathcal(E,h)_x$,
according to our assumption
$$
det^m h = \frac{1}{\prod_j |h_j|^{2n_j}} \cdot e^{-m\xi}
$$
on $\Delta$. So we can choose $g= H = \prod_j h_j^{n_j} \in H^0_{2,\det^m h}(\Delta, \det^m E)$ and
\begin{align} \label{ggg}
|g|^2_{\det^m h} = |g|^2\cdot det^m h = e^{-m\xi} \in (C_1, C_2).
\end{align}
Here $C_1, C_2$ are two strictly positive bounded real numbers. In this situation, one have
\begin{align} \label{fff}
\int_{\Delta'}|f|^2_h \cdot |g|^2_{(\det h)^m} dV_\omega  \leq  C_2 \int_{\Delta'} |f|^2_h dV_\omega < +\infty.
\end{align}
Here $\Delta'$ is the very small neighborhood of $x$. Therefore we have $f_x\cdot g_x \in \mathcal{E}(E\otimes \det^m E, h\otimes \det^m h)_x$.
According to equality \eqref{equa1}, $(e \otimes \det^m e)_x = \mathcal{E}(E\otimes \det^m E, h\otimes \det^m h)_x$, there exist a global section $r \in H^0_{2, h\otimes \det^m h}(\Delta, E\otimes \det^m E)$ such that
$$
r_x = f_x\cdot g_x.
$$

\emph{step $4$}. We choose a sufficient small neighbour $U$ of $x$ and a cut off function $\rho$ such that supp $\rho \subseteq U$ and $\rho = 1$ on a small neighbour of $x$.

Now we want to solve the $\dbar$-equation $\dbar(\rho r)= \dbar u$ with the weighted $L^2$-estimate. We first define two weighted psh functions
\begin{align*}
\varphi_k &= (n+k+k') \log |z-x|^2 + |z|^2; \\
\varphi_{k, \delta} &= (n+k+k') \log (|z-k|^2 + \delta^2) + |z|^2.
\end{align*}
Here $k'$ represent the order of $g$ at the point $x$, i.e., $g_x \in m_x^{k'}$ but $g_x \notin m_x^{k'+1}$. The previous Theorem \ref{debarrq1} shows that there exists an $E\otimes \det^m E$ valued $(0,0)$ form $u$ such that $\dbar u = \dbar(\rho r)$ and
$$
\int_{\Delta}|u|^2_{h\otimes \det^m h} e^{-\varphi_{k,\delta}} dV_\omega \leq \int_{\Delta}|\dbar(\rho r)|^2_{h\otimes \det^m h} e^{-\varphi_{k,\delta}}dV_\omega < +\infty.
$$
Taking some subsequence and taking the limit when $\delta \rightarrow 0$, we can get
$$
\int_{\Delta}|u|^2_{h\otimes \det^m h} e^{-\varphi_{k}}dV_\omega \leq \int_{\Delta}|\dbar(\rho r)|^2_{h\otimes \det^m h} e^{-\varphi_{k}}dV_\omega < +\infty.
$$
By the definition of $\varphi_k$, we know that $u_x \in m_x^{k+k'+1}\cdot (E\otimes \det^m E)_{(x)}$, this is actually not obvious. Indeed, by Lemma $2.3$ in \cite{Iwa21}, if we write $u=\sum u_i \delta_i$, here $\delta_i$ are the holomorphic frame of vector bundle $E\otimes \det^m E$, we then have each $\delta_i \in m_x^{k+k'+1}$.
Since $\dbar(\rho\cdot r - u)=0$, this means $(\rho\cdot r-u)$ is holomorphic section of $E\otimes \det^m E$. Then $(\frac{\rho\cdot r}{g}- \frac{u}{g})$ is a meromorphic section of $E$ on the whole $\Delta$ and holomorphic on $\Delta\backslash V(g)$, where $V(g):= \{z\in \Delta : g(z)=0\}$. On one hand, due to the fact
$$(\frac{\rho\cdot r}{g})_x = \frac{r_x}{g_x} = f_x$$
and $\rho$ is a cut off function near $x$ with small enough support, the first term $\frac{\rho\cdot r}{g}$ is smooth on the whole $\Delta$ and satisfying
$$
\int_{\Delta} |\frac{\rho\cdot r}{g}|^2_h dV_\omega < +\infty.
$$
On the other hand, if we have $\int_{\Delta} |\frac{u}{g}|^2_h dV_\omega < +\infty$, then $(\frac{\rho\cdot r}{g}- \frac{u}{g})$ is square integrable with respect to metric $h$, hence it can be extended to the whole $\Delta$ as a holomorphic section of $E$. Grant this for the time being, let $s= \frac{\rho\cdot r}{g}- \frac{u}{g}$, take the germ at $x$, we have $s_x = f_x - \frac{u_x}{g_x}$ and
$$
f_x = s_x + \frac{u_x}{g_x} \in e_x + \mathcal{E}(E,h)_x\cap m_x^k\cdot E_{(x)}.
$$
This is to say, $f_x \subseteq e_x + \mathcal{E}(E,h)_x \cap m_x^k\cdot E_{(x)}$ and therefore equality  \eqref{equal2} is obtained.We now prove the desired result that $\mathcal{E}(E,h)$ is coherent.

\emph{step $5$}. The last step is to show that $\int_{\Delta} |\frac{u}{g}|^2_h dV_\omega < +\infty$. But this is easy by the choose of the $g$, see \eqref{ggg}, indeed, we have
\begin{align*}
\int_{\Delta} \mid\frac{u}{g}\mid^2_h dV_\omega &= \int_{\Delta} \frac{|u|^2_{h\otimes \det^m h}}{|g|^2_{\det^m h}} dV_\omega \\
&\leq \frac{1}{C_1} \int_{\Delta} |u|^2_{h\otimes \det^m h} dV_\omega < +\infty
\end{align*}
as desired.
\end{proof}

\begin{remark} \label{3points}
Look through the proof of theorem \ref{maintheorem}, the next three conditions are sufficient:
\begin{enumerate}
\item The vector bundle $(E\otimes \det^m E, h\otimes \det^m h)$ has the $L^2$-estimate of the associated $\dbar$-equation and therefore the sheaf $\mathcal{E}(E\otimes \det^m E, h\otimes \det^m h)$ is coherent.
\item The determinant bundle $(\det^m E, \det^m h)$ is singular Griffiths semi-positive when $(E, h)$ is singular Griffiths semi-positive. Thus $\mathcal{E}(\det^m E, \det^m h)$ is coherent, this is a classical result. Note $m$ can be any positive integer.
\item If the weight of metric $\det h$ has algebraic singularities, then we choose some integer $m$ such that $-m \log \det h = \sum_j n_j \log |h_j|^2 + m\xi$ just like \eqref{algebraic}.
\end{enumerate}

\end{remark}

\section{Singular Hermitian metrics with analytic singularities}

We just prove the coherence of $\mathcal{E}(E,h)$ when the weight $\psi=-\log \det h$ has algebraic singularities. It seems the conclusion should be valid for the case that $\psi$ has analytic singularities. We shall prove it now.

\begin{theorem} \label{debarrq2}
Let $(\Delta, \omega)$ be a small polydisk with a standard metric $\omega$, and $\phi$ be a smooth strictly plurisubharmonic function on $\Delta$.
let $(E,h)$ be the trivial holomorphic vector bundle with a Griffiths semi-positive singular Hermitian metric $h$. There exists a natural induced metric $h\otimes (\det h)^m$ on the vector bundle $E\otimes \det^m E$ for any $m \in \mathbb{N}$.  If there is a psh function $\theta$ such that $\xu\pa\dbar \theta - m\xu\pa\dbar(-\log(\det h)) \geq 0$, then for any $\dbar$-closed $E\otimes \det^m E$ valued $(n,q)$-form $u$ satisfying
$$
 \int_{\Delta} |u|^2_{\omega, h\otimes e^{-\theta}} \cdot e^{-\phi} dV_{\omega} \leq +\infty.
$$
We can find a $E\otimes \det^m E$ valued $(n,q-1)$-form $\alpha$ such that $\dbar \alpha = u$ and
$$
\int_{\Delta} |\alpha|^2_{\omega, h\otimes e^{-\theta}}\cdot e^{-\phi} dV_{\omega} \leq C \int_{\Delta} |u|^2_{\omega, h\otimes e^{-\theta}} e^{-\phi} dV_{\omega}
$$
for some constant $C$.
\end{theorem}

\begin{remark}
Here we make $e^{-\theta}$ the non-trivial metric on the trivial line bundle $\det^m E$. This theorem tells us that the vector bundle $(E\otimes det^m E, h\otimes e^{-\theta})$ has the $L^2$-estimate of the associated $\dbar$-equation and therefore the sheaf $\mathcal{E}(E\otimes det^m E, h\otimes e^{-\theta})$ is coherent.
\end{remark}

\begin{proof} [The proof of Theorem \ref{debarrq2}]
Like the proof in the previous Theorem \ref{debarrq1}, we first regularize the metric $h$ and the psh function $\theta$ on the smaller polydisk as Lemma \ref{regular}. The approximation sequence of smooth metrics $h_\nu$ increasing to $h$ and the approximation sequence of smooth psh functions $\theta_\nu$ decreasing to $\theta$.
We then calculate the curvature of trivial vector bundle $(E\otimes \det^m E, h_\nu\otimes e^{-\theta_\nu}\cdot e^{-\phi})$. Firstly we write $h_\nu \otimes e^{-\theta_\nu} = h_\nu\otimes \det^m h_\nu \cdot ( \det^m h_\nu)^{-1}\cdot e^{-\theta_\nu}$.

%From the result of BP and Raufi, we can obtain a sequence of smooth Hermitian metrics $\{h_v\}_{v=1}^{\infty}$ with Griffiths semi-positive curvature increasing to $h$ on any relatively compact subset in $\Delta$.

Due to Demailly-Skoda theorem, $h_\nu \otimes \det^m h_\nu$ is Nakano semi-positive. At the same time, we can arrange the sequence smooth psh function $\theta_\nu$ such that
$$\xu\pa\dbar \theta_\nu - m\xu\pa\dbar(-\log(\det h_\nu)) \geq -C \omega,
$$
here $C$ are some positive constant independent to $\nu$. Hence the curvature of vector bundle $(E\otimes \det^m E, h_\nu\otimes e^{-\theta_\nu} \cdot e^{-\phi})$ can be calculated as follow
\begin{align*}
\xu \Theta_{h_\nu \otimes e^{-\theta_\nu} \cdot e^{-\phi}} = &\xu \Theta_{h_\nu\otimes \det^m h_\nu} + \xu \pa\dbar (\theta_{\nu}+ \log (\det h_\nu)^m) \otimes \Id_{E\otimes \det^m E} \\
&+\xu \pa\dbar\phi \otimes \Id_{E\otimes \det^m E}.
\end{align*}
Now if we choose the very positive psh function $\phi$, we can make the curvature of vector bundle $(E\otimes \det^m E, h_\nu \otimes e^{-\theta_\nu} \cdot e^{-\phi})$ is strictly Nakano positive, i.e.,
\begin{align*}
 &\xu \Theta_{h_\nu\otimes \det^m h_\nu} + \xu \pa\dbar (\theta_{\nu}+ \log (\det h_\nu)^m) \otimes \Id_{E\otimes \det^m E} +\xu \pa\dbar\phi \otimes \Id_{E\otimes \det^m E} \\
&\geq \varepsilon \omega\otimes Id_{E\otimes \det^m E}
\end{align*}
Thus for any $E\otimes \det^m E$ valued $(n,q)$-form $u$, we have
$$
\langle  [\xu \Theta_{h_\nu \otimes e^{-\theta_\nu} \cdot e^{-\phi}}, \Lambda_{\omega}]u, u \rangle \geq q\varepsilon |u|^2.
$$
The rest proof basically the same as the previous Theorem \ref{debarrq1} and we omit it.
\end{proof}

\begin{corollary} \label{theta}
With the same condition in the Theorem \ref{debarrq2}, the sheaf $\mathcal{E}(E\otimes \det^m E, h\otimes e^{-\theta})$ is coherent.
\end{corollary}

\begin{remark} \label{4points}
Similar to Remark \ref{3points}, once the next three conditions are satisfied,
then from the proof of Theorem \ref{maintheorem}, we know $\mathcal{E}(E, h)$ is coherent.
\begin{enumerate}
\item The vector bundle $(E\otimes \det^m E, h\otimes e^{-\theta})$ has the $L^2$-estimate of the associated $\dbar$-equation and therefore the sheaf $\mathcal{E}(E\otimes \det^m E, h\otimes e^{-\theta})$ is coherent.
\item The determinant bundle $(\det^m E, e^{-\theta})$ is singular Griffiths semi-positive. Thus $\mathcal{E}(\det^m E, e^{-\theta})$ is coherent.
\item If the weight $\theta$ has algebraic singularities, and moreover the coefficient is positive integer, i.e., $\theta= \alpha \log (|f_1|^2 + \cdot\cdot\cdot + |f_N|^2) + v$, for $\alpha \in \mathbb{Z}^+$ and holomorphic functions $f_i$ and bounded function $v$.
\end{enumerate}
\end{remark}

\begin{corollary} \label{maintheorem1}
Let $(E,h)$ be a holomorphic vector bundle over an $n$-dimensional complex manifold with a singular Griffiths semi-positive Hermitian metric $h$. If the weight of metric $\det h$ has analytic singularities, then $\mathcal{E}(E,h)$ is coherent.
\end{corollary}

\begin{proof}
Locally $(-\log\det h)$ can be written as
$$
-\log\det h = \alpha \log (|f_1|^2 + \cdot\cdot\cdot + |f_N|^2) + v
$$
where $\alpha$ is a positive real number, $v$ is a locally bounded function and all $f_i$ are holomorphic functions. We can choose a positive integer $\beta$ such that $\beta \geq \alpha$. It is obviously that $\theta := \beta \log (|f_1|^2 + \cdot\cdot\cdot + |f_N|^2) + v$ is plurisubharmonic function and satisfying $\xu \pa\dbar \theta - \xu\pa\dbar (-\log \det h) \geq 0$. By Theorem \ref{debarrq2} and Corollary \ref{theta}, the sheaf $\mathcal{E}(E\otimes \det E, h\otimes e^{-\beta})$ is coherent. According to the above Remark \ref{4points}, we can obtain the desired result, i.e., the sheaf $\mathcal{E}(E,h)$ is coherent.
\end{proof}

\begin{remark}
It is natural to try to prove the general case by a certain kind of approximation, but there is one difficulty.
I want to approximate $ \det h $ with a sequence of psh functions $ f_i$ with analytic singularities. Moreover according to Theorem \ref{debarrq2}, one may need the next result:
$ \xu \pa\dbar f_i - \xu \dbar \det h $ are uniformed bounded below, i.e., $ \xu \pa\dbar f_i - \xu \pa\dbar \det h \geq -C \omega $, where $C$ is a constant independent to $i$.
Only by Demailly's approximation, I think this is hard to make it. Demailly's approximation shows that if
  $\xu \pa\dbar \det h \geq \gamma $  for some smooth $(1,1)$-forms $\gamma$, then we have $ \xu \pa\dbar f_i - \gamma \geq -C \omega $  for some constant $C$.
This is not enough for our purpose.
%So far, I don't know how to solve this problem.
\end{remark}

\end{document}